\documentclass[12pt]{article}
\usepackage[cp1251]{inputenc}
\usepackage[english]{babel}
\usepackage[final]{epsfig}
\usepackage{amsfonts, epsfig, float}
\usepackage{graphicx}
\usepackage{amssymb,latexsym}
\usepackage{amsmath}
\usepackage{stmaryrd}
\usepackage{textcomp}
\usepackage{epsfig, float}

\usepackage{color}
\newtheorem{proposition}{Proposition}
\newtheorem{theorem}{Theorem}
\newtheorem{remark}{Remark}

\newtheorem{lemma}{Lemma}
\newtheorem{definition}{Definition}

\newtheorem{corollary}{Corollary}
\newcommand{\be}{\begin{equation}}
\newcommand{\ee}{\end{equation}}
\newcommand{\bea}{\begin{eqnarray}}
\newcommand{\eea}{\end{eqnarray}}

\newcommand{\cop}{{\cal COP}^p}
\newcommand{\ck}{{\cal K}}
\newcommand{\cd}{{\cal D}}

\begin{document}

%\articletype{Research Paper}% Specify the article type or omit as appropriate

\title{On equivalent representations and  properties of faces of the cone of copositive matrices }

\author{Kostyukova O.I.\thanks{Institute of Mathematics, National Academy of Sciences of Belarus, Surganov str. 11, 220072, Minsk,
 Belarus  ({\tt kostyukova@im.bas-net.by}).} \and Tchemisova T.V.\thanks{Mathematical Department, University of Aveiro, Campus
Universitario Santiago, 3810-193, Aveiro, Portugal ({\tt tatiana@ua.pt}).}}

\maketitle

\begin{abstract}
The paper is devoted to a study of the cone $\cop$  of copositive matrices. Based on the known from semi-infinite optimization concept of immobile indices, we define zero and minimal zero vectors of a subset of the cone $\cop$  and use them to obtain different  representations of faces of $\cop$  and the corresponding dual cones. We describe the minimal face of  $\cop$  containing a given convex subset of this  cone and prove some propositions  that allow to obtain equivalent   descriptions of the feasible sets of a copositive problems and may be useful for creating new numerical methods based on their regularization.
\end{abstract}

\\


\section{Introduction}\label{Introduction. Motivation}

This work is motivated by our main challenge: study of  Copositive Programming (CoP) problems and their properties.  CoP deals with  a special class of conic  problems and can be considered as an optimization over the convex cone of so-called {\it copositive matrices} (i.e. matrices which are positive semidefinite on the non-negative orthant). Copositive problems attract the attention of researchers as they have many interesting uses (see e.g.  \cite{AH2013,Bomze2012,Dur2010}, and the references therein). It should be noted here that CoP can be seen as a generalization of Semidefinite Programming (SDP) and  a  particular case   of Semi-infinite  Programming (SIP), whose important applications are well known \cite{HandbookSDPnew, Weber2, W-handbook}.

An optimization problem is called {\it regular} if its constraints satisfy some additional (regularity) conditions, so-called {constraint qualifications (CQs).} The regularity of an optimization problem  guarantees that the first-order necessary optimality conditions can be formulated in the  form of so-called Karush-Kuhn-Tucker (KKT) system \cite{Liu, Lop-Still}  and the strong duality relation is satisfied \cite{Kortanek}. In most of the works on CoP {and SIP}, the study is conducted on the assumption that some regularity conditions are satisfied (see e.g. \cite{AH2013}  for linear CoP and \cite{M-F, Klatte2,  Morduck} for convex SIP).

In the cases where  regularity conditions are not met, there is no guarantee that  an  optimal solution satisfies the KKT  type optimality conditions, the  first-order optimality conditions of the Fritz John type (see Theorem 2 in \cite{Lop-Still}) become   degenerate (they are satisfied for all feasible solutions and hence are not informative),  and the strong duality relation may fail. This creates difficulties in numerical solving problems. To overcome these obstacles, special regularization procedures can be applied, allowing to rewrite the original    problem in an equivalent form for which certain regularity conditions are met. This approach is based on the fact that as a rule, the violation of CQs is caused by an unhappy description of the feasible set of optimization problem.
  To obtain non-degenerate optimality conditions, one needs a "good" representation of the feasible set.

 In conic optimization, the regularization procedures are based on the so-called {\it facial reduction algorithms} (FRAs). To perform a constructive regularization procedure  for a conic  problem, it is  necessary to  explicitly describe certain faces of its feasible cone (the faces containing a given convex subset of this cone) and their dual ones.
For SDP problems, which are a particular case of conic problems, the properties of faces of the cone of semidefinite matrices are well studied. This
allowed to describe  constructive regularization procedures for SDP and develop a duality theory  satisfying strong duality conditions without regularity assumption. In particular, in \cite{Ramana, Ramana-W}, etc.,   a  dual problem for SDP is  explicitly formulated in the form of an {\it Extended Lagrange Dual Problem} (ELDP). Several attempts have been done to obtain regularization procedures for general conic problems (see e.g. \cite{Dry-Wol, Waki}), but the  procedures  described there  are implicit  and do not permit to obtain explicit strong dual formulations. There is no such duality theory for CoP as well. This can be explained by the fact that the structure and properties of faces of the  cone of  copositive matrices  and the respective duals are not well studied yet \cite{Dur}. It worth mentioning that several attempts have been  done to study the facial structure of the cones of copositive and completely positive matrices. Thus in  \cite{Afonin, Dick, Hild},  the authors give explicit characterizations of extreme rays (faces of dimension one) of  copositive cones  of   dimensions  five and six. In \cite{Dick1, Dick},   some properties of special types of faces (minimal and maximal ones) are studied. Nevertheless,    until now, for the cone of copositive matrices,    faces  of this cone  and  the corresponding  dual cones are not well studied.

%\vspace{3mm}

 The aims of the paper are as follows:

\vspace{1mm}

$\bullet$ for a given convex subset of the set of copositive  matrices, to define the corresponding set of all its zeros and minimal zeroes    and to study some of their properties;

 $\bullet$ given a face  of the cone of copositive matrices, to obtain its  representation  in terms of the corresponding minimal  zeros  and  describe explicitly the dual cone to this face;

 $\bullet$ for a given  closed convex  subset of the cone of copositive matrices, to derive an  explicit description of the  cone's  minimal face  containing this subset  in terms of the corresponding set of  minimal   zeros;

$\bullet$  to  obtain equivalent ``good'' descriptions of the feasible sets of  convex CoP problems.

\vspace{2mm}

The rest of the paper is structured as follows. In section \ref{S2}, we collect notation,  basic definitions, and prove some important results which will be used in this paper. In section \ref{S4}, for   a convex closed subset of   the cone of copositive matrices, we define  its {\it minimally active}  element and prove the existence of  such an element. Section \ref{S5} contains three equivalent  representations of  faces    of the copositive cone, and in section \ref{S6}, we deduce alternative representations of the dual cones to the faces. In section \ref{S7}, we describe the minimal face of the cone of copositive matrices containing a given convex set  and prove some corollaries  that may be useful for creating new numerical methods based on the minimal cone representations and  regularization procedures for  copositive problems. In section \ref{S3}, we  consider equivalent   descriptions of the feasible set of a copositive problem. The final section \ref{S8} contains some conclusions.

\section{ Notation, basic definitions, and  important preliminary results}\label{S2}

Given an integer $p>1$, consider  the vector space $\mathbb R^p$ with the standard orthogonal basis   $\{e_k,  \ k =1,2, \dots, p\}$. Denote by  $\mathbb R^{ p}_+$  the set of all $ p$ - vectors with non-negative components,
by ${\mathcal S}(p)$ and $\mathcal S_+(p)$  the space of real symmetric $p\times p$ matrices
and the cone of symmetric positive semidefinite $p\times p$ matrices, respectively, and let
$\mathcal{COP}^{p}$  stay for the cone of symmetric copositive $p\times p$ matrices:
$$\mathcal{COP}^p:=\{D\in {\mathcal S}(p):t^{\top}Dt\geq 0\ \forall t \in \mathbb
R^p_+\}.$$ The space $\mathcal S(p)$ is considered here as a vector space with the trace inner product $A\bullet B:={\rm trace}\, (AB).$
 Denote \begin{equation}\label{SetT}T:=\{t \in \mathbb R^p_+:\mathbf{e}^{\top}t=1\} \end{equation}
with $\mathbf{e}=(1,1,...,1)^{\top}\in \mathbb R^p $. It is evident that the cone $\cop$ can be equivalently described as follows:
$\cop=\{D\in \mathcal S(p): \ t^{\top}Dt\geq 0 \ \forall t\in T\}.$

Given a vector  $t=(t_k,k\in P)^{\top}\in \mathbb R^p_+$ with   $P:=\{1,2,...,p\}$,  introduce the  sets
$$P_+(t):=\{k\in P:t_k>0\}\  \mbox{ and } \ P_0(t):=P\setminus P_+(t).$$
Given a set $\mathcal B\subset\mathbb R^p$  denote by
${\rm conv } \mathcal B$   the convex hull of  $\mathcal B$,  by  ${\rm int  }\,\mathcal B$   the interior, by ${\rm relint  }\,\mathcal B$ the relative  interior and by ${\rm cl  }\,\mathcal B$ the closure  of this set.
For the set $\mathcal B$  and a point $l=(l_k,k \in P)$,  denote by $\rho(l,\mathcal B)$ the distance between this set and the point,
$\ \rho(l,\mathcal B):=\min\limits_{\tau \in \mathcal B}\sum\limits_{k\in P}|l_k-\tau_k|.$

In our previous papers (see e.g. \cite{KT-JOTA,KT-SetValued,KT-SOIC}, and the references therein), considering  convex  SIP problems, for a constraint function satisfying the inequalities $f(x,t) \geq 0$ $  \forall x\in X_f, $ $  \forall  t\in \Upsilon,$
where $X_f\subset \mathbb R^n$ is the convex set of feasible solutions, $\Upsilon\subset \mathbb R^p$ is a compact index set, we introduced the concept of immobile indices of the constraints
and  showed that this  concept plays an important role in the study of optimality  of  solutions to these problems and regularity of their feasible sets.
The set of immobile indices was defined as follows:
$$T_{im}:=\{t \in \Upsilon:f(x,t)=0\;\forall x \in X_{f}\}.$$
In this paper, in a similar way we introduce a notion of  zeros of convex subsets of the cone $\cop$ and use them to describe faces   of this cone, the dual cones corresponding to faces of $\cop$, and the minimal face containing a given convex subset.

For a matrix set $Q \subset \cop$ we will say that vector $t\in T$ is a {\it zero  of the set} $Q$ if
$t^{\top} Dt=0 \  \forall D\in Q.$

Denote by $T_0(Q)$ the {\it set of all zeros} of the given matrix set $Q$:
\begin{equation}\label{Formula2} {T_0}(Q):=\{t \in   T : t^{\top} Dt=0 \ \  \forall D \in Q\}.\end{equation}

Note that if  we consider a function $\varphi(D,t):=t^\top Dt$ defined for $ D\in {\cal S}(p), t \in T$,  and satisfying the following  conditions:   $$\varphi(D,t) \geq 0\  \forall   D \in Q, \ \forall t\in T,$$ where $Q \subset \cop $,   then  $t\in T $ is an {\it immobile index} of this function if $t^{\top}Dt=0 \ \forall D\in Q.$
Therefore, the set ${T_0}(Q)$ of all zeros  of the matrix set $Q$ is nothing but the set of immobile indices of the defined above function $\varphi$.

 For a given subset $Q \subset \cop $, a vector $ \tau\in  T_0(Q)$ is called a {\it minimal zero} of $Q$ if there does not exist $t \in T_0(Q)$, $t\not =\tau,$  such that $P_+(t)\subset P_+(\tau).$

Note that  the   introduced here definitions  of zeros and minimal zeros of a given set $Q \subset \cop,$  generalize the concepts of zeros and minimal zeros introduced in \cite{Dick}  for a given matrix $A\in \cop.$

For the considered above set $Q$, let us say that it satisfies {\it the  Slater condition} if
\begin{equation}\label{Slait1}
\exists \; \bar{D} \in Q\; \mbox{ such that } \;\bar D \in {\rm int }\, \mathcal{COP}^p=
\{D\in \mathcal{S}(p): \  t^{\top}Dt>0\ \forall t \in T\}.
\end{equation}

Repeating the chain of proofs of  Lemma 1 and  Proposition 1  in \cite{KT-SetValued}, we can prove the following lemma.

\begin{lemma} \label{LOK} Given a closed convex subset $Q\subset \cop $,
\begin{enumerate}
  \item[(i)] the Slater condition (\ref{Slait1}) is equivalent to the emptiness of the  set $ T_0(Q)$;
 \item[(ii)] the  set $ T_0(Q)$ of all zeros of the matrix set $Q$  is either empty or  can  be represented as a union of a finite number of convex  closed  bounded polyhedra.
    \end{enumerate}
\end{lemma}

It follows from Lemma \ref{LOK} that if  $ T_0(Q)\neq \emptyset$, then the set ${\rm conv}\, T_0(Q) $ is a convex bounded polyhedron with a finite number of vertices.  If $T_0(Q)=\emptyset$, then, evidently,   ${\rm conv} T_0(Q)=\emptyset$.
 Denote by
\be Z_Q:=\{\tau(j), j \in J_Q\},\; 0\leq |J_Q|<\infty, \label{vert}\ee
the set of all vertices of the set ${\rm conv}\, {T_0}(Q) $. Notice that $J_Q:=\emptyset$ when ${T_0}(Q)=\emptyset$.

From the definition of the set ${T_0}(Q)$, we conclude that for each matrix $D\in Q$, the {elements}   $\tau(j) , i\in J_Q,$   of the set $Z_Q$ are optimal solutions of the problem
\be \min \; t^\top Dt\;\; {\rm s.t. } \; t \in T.\label{1problem}\ee

\vspace{-3mm}

Hence
\be e^\top_kD\tau(j)=0 \ \forall k \in P_+(\tau(j)); \;  e^\top_kD\tau(j)\geq 0 \ \forall  k \in P_0(\tau(j)), \  \forall j \in J_Q, \; \forall D \in Q.\label{1in}\ee

 Define the sets
\be M_Q(j):=\{k\in P: e^\top_kD\tau(j)=0  \ \forall D\in Q\}, \; j \in J_Q.\label{M}\ee
It follows from (\ref{1in}) that $P_+(\tau(j))\subset M_Q(j)\ \forall j \in J_Q.$

Note that in this paper we consider different types of subsets $Q$ of the cone $\cop$. Therefore we include the reference on $Q$ in the designation of the sets (\ref{Formula2}), (\ref{vert}), and (\ref{M}).

\begin{lemma}\label{Lmin11}
For a given convex closed set $Q\subset\cop$, the set of all its minimal zeros  is either empty or  finite and coincident with the set $\{\tau(j), j \in J_Q\}$  of all vertices of the set ${\rm conv}\, T_0(Q) $.
\end{lemma}

{\bf Proof}. It is evident that if the set $T_0(Q)$ is empty, then  the set of minimal zeros is empty too.

 Suppose that $T_0(Q)\not =\emptyset$. Notice that since $Q\subset {\cal COP}^p$ and $\tau(j)\in T_0(Q)$ for all $ j \in J_Q$,  the  relations (\ref{1in}) hold true.
%\be e^\top_kA\xi(j)=0, \; k \in P_+(\xi(j)); \; e^\top_kA\xi(j)\geq 0, \; k \in P \setminus P_+(\xi(j))\;\; \forall  j \in J.\label{NN3}\ee

First, let us show that all  vertices of the set ${\rm conv } T_0(Q)$   are  minimal  zeros of the set $Q$.
Suppose that, on the contrary, for some $j_0\in J_Q$, {the corresponding} zero $\tau(j_0)$ is not minimal for the set $Q$.  Then there exists another zero $\bar t \in T_0(Q)$ such that $P_+(\bar t)\subset P_+(\tau(j_0))$, $\bar t \not=\tau(j_0).$
 Consequently, for a sufficiently small $\theta>0$,  we have
 \be \bar \tau :=(\tau(j_0)-\theta \bar t)/(1-\theta)\in T, \;\; \bar \tau \not =\bar t.\label{bar}\ee
 Let us show that $\bar t^{\, \top} D\tau(j_0)=0$ for all $D\in Q.$ In fact, taking into account (\ref{1in}), we get
 $$\bar t^{\, \top}D\tau(j_0)=\sum\limits_{k\in P}\bar t_ke^\top_k D\tau(j_0)=\sum\limits_{k\in P_+(\tau(j_0))}\bar t_ke^\top_k D\tau(j_0)=0\;\; \forall D\in Q.$$
For any $D\in Q$,  from (\ref{bar}) and the equalities $\bar t^{\, \top} D\bar t=(\tau(j_0))^\top D\tau(j_0)=\bar t^{\, \top} D\tau(j_0)=0$ it follows  that $\bar \tau^\top D\bar \tau=0$.  Hence
$\bar \tau \in T_0(Q). $ As a result, we obtain
$$\tau(j_0)=(1-\theta) \bar \tau +\theta \bar t, \; \theta \in (0,1), \; \bar \tau\in {\rm conv} T_0(Q), \; \bar t\in {\rm conv} T_0(Q),\; \bar t\not =\bar \tau.$$
But these relations contradict the assumption that $\tau(j_0)$ is a vertex of the set  ${\rm conv} T_0(Q)$.
Thus, it is proved  that  all vertices $\tau(j), \ j\in J_Q,$ are  minimal zeros of the set  $Q$.

\vspace{1mm}

 To  prove that  there isn't a single minimal zero outside  the  set $\{\tau(j), \ j\in J_Q\}$, c onsider a zero
$t^*\in T_0(Q)\setminus\{\tau(j), j \in J_Q\}\subset {\rm conv} T_0(Q) .$  By construction, it admits  the following representation:
$$t^*=\sum\limits_{j\in J^*}\alpha_j\tau(j),\;   \sum\limits_{j\in J^*}\alpha_j=1,\; \alpha_j>0\ \forall j \in J^*\subset J_Q,\;|J^*|\geq 2.$$
It follows from this representation that $P_+(t^*)=\bigcup\limits_{j \in J^*}P_+(\tau(j))$, which implies
$$P_+(\tau(j))\subset P_+(t^*),  \; \tau(j)\not= t^*,\;\; \tau(j)\in T_0(Q)\;\;  \forall  j \in J^*.$$
By definition, the obtained inclusions mean  that $t^*$ is not a minimal zero of the set $Q$.
The proposition is proved. $ \;\blacksquare$

\vspace{3mm}

In what follows, the set $Z_Q$ will be also called the {\it set of minimal zeros} of the set $Q.$

\vspace{3mm}

Consider the  set $T$  defined in (\ref{SetT}).
For a  given  nonempty finite subset
\be V:=\{t(i)\in T, \; i \in I\},\; 0<|I|<\infty,\label{0**}\ee
 define the number

 \vspace{-7mm}

\begin{equation} \sigma (V):=\min\{t_k(i),\;  k \in P_+(t(i)), \; i \in I\}.\label{o-1}\end{equation}
 By definition, evidently,  $\sigma (V)>0$.

Introduce the following sets:
\bea& {\Omega}(V):=\{t \in T: \rho(t,{\rm conv} V)\geq \sigma (V)\},\label{set-hat1}\\
& \mathcal N(V):=\{t \in T:\rho(t,{\rm conv }\, V)\leq \sigma (V)\}.\label{o-2}\eea
Note that, by construction, $T={\Omega}(V)\cup \mathcal N(V)$ and ${\Omega}(V)\cap \mathcal N(V)=\{t \in T: \rho(t,{\rm conv} V)=\sigma (V)\}$.

In the  rest of this section, we will prove two auxiliary  statements that will be used  to justify the main results of this paper.
\begin{proposition}\label{Pnew0}
Let the set  $V$ be given in (\ref{0**}), and the corresponding number  $\sigma (V)$ and  set $\mathcal N(V)$  defined in  (\ref{o-1}) and
 (\ref{o-2}).
 Then for any $t\in \mathcal N(V)$, there exists  $t(i_0)\in V$  such that $P_0(t)\subset P_0(t(i_0)).$
\end{proposition}
{\bf Proof.} Consider $t \in  \mathcal N(V)$. By construction, there exists $\tau \in {\rm conv}\, V$ and a nonempty set $I_*\subset I$ such that
\be \rho(t,{\rm conv} V)=\rho(t,\tau)\leq \sigma(V);\;\;
 \tau=\sum\limits_{i\in I_*}\alpha_it(i), \; \sum\limits_{i\in  I_*}\alpha_i=1, \; \alpha_i>0 \; \forall i \in  I_*.\label{o-5}\ee

Suppose that $P_0(t)\not \subset P_0(t(i)) \; \; \forall i \in I. $  It is evident that these relations are equivalent to the following ones: $P_0(t)\cap P_+(t(i))\not= \emptyset \; \; \forall i \in I. $
Taking into account the latest inequalities and (\ref{o-5}),  we get
\begin{equation*}\begin{split}\sigma (V)&\geq \rho(t,\tau)=\sum\limits_{k\in P}|t_k-\tau_k|=\!\!\sum\limits_{k\in P_+(t)}|t_k-\tau_k|+\sum\limits_{k\in P_0(t)}\tau_k\\
&\geq\sum\limits_{k\in P_0(t)}\tau_k= \sum\limits_{k\in P_0(t)}\sum\limits_{i\in  I_*}\alpha_it_k(i)=\sum\limits_{i\in  I_*}\alpha_i\sum\limits_{k\in P_0(t)\cap P_+(t(i))}t_k(i)\geq \sigma (V).\end{split}\end{equation*}
Thence,  $\rho(t,\tau)=\sigma (V)$ and in the chain of inequalities above, we can replace the inequality symbol  by the equality one. Then
$\sum\limits_{k\in P_+(t)}|t_k-\tau_k|=0$. Hence $t_k=\tau_k  $ $ \forall {k\in P_+(t)}.$
From the latter equalities and the following relations: $$ t \geq 0, \ \tau\geq 0,\ {\mathbf e}^{\top}t={\mathbf e}^{\top}\tau=1; \  t_k=0\; \ \forall k  \in  P_0(t),$$ we get the equality $t=\tau$. Then, taking into account relations (\ref{o-5}), we conclude that $P_0(t)= \bigcap\limits_{i\in I_*}P_0(t(i))$, which
implies $P_0(t)\subset P_0(t(i))\; \forall i \in  I_*.$

Hence we obtain a contradiction  with the assumption $P_0(t)\not \subset P_0(t(i))\    \forall i \in I.$
The proposition is proved. $ \blacksquare$

\begin{theorem}\label{LO-1} Consider the    defined in (\ref{SetT})   set $T$  , any its subset  $V$  in the form (\ref{0**}), and the corresponding set $\Omega(V) $ defined in (\ref{set-hat1}).
 For any matrix $D\in S(p)$, the relations
\begin{equation} t^{\top}Dt\geq 0\; \forall  t \in \Omega(V){\mbox{ and }} \;  Dt(i)\geq 0 \; \forall i \in {I},\label{4*}\end{equation}
%&e^T_kD\tau(i)=0, k \in P_+(\tau(i)),\; e^T_kD\tau(i)\geq 0, k \in P_0(\tau(i)),\; i \in I,\label{4*-1}\eea
imply
\be t^{\top}Dt\geq 0 \; \forall  t \in T.\label{5*}\ee
\end{theorem}

{\bf Proof.} Given a matrix $D\in S(p)$ and the set $V$ defined in (\ref{0**}), first notice that  the relations $ Dt(i)\geq 0 \; \forall i \in I,$ imply
$ t^{\top}Dt\geq 0\;\; \forall  t \in {\rm conv} \{t(i), i \in I\}={\rm conv}V.$
%\label{4*-2}\ee

Suppose that  for a given $D\in S(p)$, relations (\ref{4*}) hold true, but some of the inequalities in  (\ref{5*}) are violated.  Hence for a vector  $\bar t$ given by
\begin{equation} \bar t: ={\rm arg} \{\min t^\top Dt,\; \mbox{ s.t. } \; t \in T\},\label{6*}\end{equation}
we have
$\bar t^{\top}D\bar t<0 \mbox{ and } \bar t\in \mathcal N(V)\setminus{\rm conv} V,$
where the set $\mathcal N(V)$ {is} defined in (\ref{o-1}), (\ref{o-2}). Since $\bar t \in \mathcal N(V), $  from Proposition  \ref{Pnew0} it follows:
\be \exists  \; i_0\in I \mbox{  such that } P_0(\bar t)\subset P_0(t(i_0)).\label{7*}\ee
 Let us set $l:=\bar t -t(i_0)$. The vector $l$ is  a feasible  direction for $t(i_0)$ in the set $T$ as, evidently,
$t(i_0)+\lambda(\bar t -t(i_0))\in T$ for all $\lambda\in [0,1].$
It follows from (\ref{7*}) that $ l_k=0 \; \; \forall  k \in P_0(\bar t).$ Hence, the vector $l$ is  a feasible direction for $\bar t$ in $T$ as well. As $l$ is  a feasible  {direction} for $t(i_0)$  and $\bar t$  in $T$, then there exists $\gamma_0>1$  such that for all $\gamma\in [0,\gamma_0]$, it holds:
$t(\gamma):=t(i_0)+\gamma l=t(i_0)+\gamma(\bar t -t(i_0))\in T.$

Define the function $$f(\gamma):=t^{\top}(\gamma)Dt(\gamma)=\gamma^2a+2\gamma b+c, \; \gamma\in [0,\gamma_0],$$
where $a{:}=l^{\top}Dl,$ $b{:}=l^{\top}Dt(i_0), $ and $c{:}=(t(i_0))^{\top}Dt(i_0).$
According to (\ref{6*}), we have
$$\min\limits_{0\leq \gamma\leq \gamma_0}{f(\gamma)}=f(\gamma^*=1)=a+2b+c.$$
Since $\gamma^*=1\in (0,\gamma_0)$, then  $0=\frac{df(\gamma^*)}{d\gamma}=2a+2b$, which implies $a=-b$  and, equivalently, $ l^{\top}Dl=-l^{\top}D t(i_0)$. Then,  taking into account the definition of vector $l$, we conclude that
\be \bar t^{\top}D\bar t=\bar t^{\top}Dt(i_0).\label{7**}\ee
Remind that, by assumption, it holds $\bar t^{\top}D\bar t<0$. On another hand, the inequalities $Dt(i_0)\geq 0 $ and $ \bar t\geq 0$ imply  $ \bar tDt(i_0)\geq 0.$ The obtained contradiction  with the equality  (\ref{7**}) completes the proof. $ \qquad \blacksquare$

\section{Minimally active elements of a set $Q\subset \cop$}\label{S4}

Given a matrix $A\in \cop$, let $T_0(A)$ be the  set of its zeros defined in (\ref{Formula2}) with $Q$ replaced by $A$.
Denote by
$$Z_A:=\{\tau(j)\in T_0(A), \; j \in J_A)\},  \ |J_A|<\infty,$$
the set of all  vertices (extremal points)  of the set ${\rm conv}\,T_0(A)$  which coincides with the set of minimal zeros of  $A$ (see Lemma \ref{Lmin11}), and introduce the sets
$$M_A(j):=\{k \in P:e^\top_kA\tau(j)=0\},\; j \in J_A.$$

%Notice that  $P_+(\tau(i,D))\subset L_a(i,D),\; i \in I(D).$
Let $Q$ be a  convex closed subset of the cone ${\cal COP}^p$.
\begin{definition}\label{MinAct}
A matrix $A\in Q$ is called minimally active element of the set $Q $ if  for any $D \in Q $, it holds $$ T_0(A)\subset  T_0(D) \mbox{ and } T_0(A)= T_0(D) \;  \Longrightarrow\;  M_A(j)\subset M_D(j) \  \forall j\in J_A=J_D.$$
\end{definition}
Note that from the  definition, it follows that  if $A\in Q$ and $T_0(A)=\emptyset$, then $A$ is a  minimally active element of the set $Q$.

The main result of this section  is the  proof  of  a theorem  which ensures the existence of a minimally active element of any convex closed subset of the cone ${\cal COP}^p$.

\vspace{2mm}

 Given a convex closed  matrix set $Q\subset \cop$, consider the sets   $ T_0(Q)$  and $ Z_Q $  of all zeros and minimal zeros of $Q$, and  the sets
$  M_Q(j), j \in J_Q,$ defined in (\ref{Formula2}), (\ref{vert}), and (\ref{M}).
\begin{theorem}\label{L25-2-2} Given a convex closed  set $Q\subset \cop$, suppose that the set  $T_0(Q)$  of all zeros of $Q$ is not empty.
	Then there exists a matrix    $\bar{D} \in Q$ such that
   \begin{equation}\label{For2}\begin{split}& t^{\top}\bar D t>0 \ \forall t\in T\setminus T_0(Q); \ t^{\top}\bar D t=0 \ \forall t\in T_0(Q);\\
   &e^{\top}_k\bar D\tau(j)>0 \ {\forall} k\in P\setminus M_Q(j), \ \forall j\in J_Q. \end{split}\end{equation}
	\end{theorem}
{\bf Proof.}
Given a convex closed set $Q\subset \cop$,  consider  the set $Z_Q=\{\tau(j),j\in J_Q \}$  of all its minimal zeros.
 Denote
\begin{equation}\label{N*}\begin{split}
&U:=\{(i,j): \ i\in J_{Q},\;  j\in J_{Q}, \; i\neq j\},\\
& U_{+}:=\{(i,j)\in U: \exists \; D(i,j)\in Q,\;  (\tau(i))^{\top}D(i,j)\tau(j)>0\}.
\end{split}\end{equation}
By definitions (\ref{M}) and (\ref{N*}), for any $j \in J_Q$ and any $k\in P\setminus M_Q(j)$ there exists a matrix $A(k,j)\in Q$ such that $e^{\top}_kA(k,j)\tau(j)>0,$
and for any $(i,j)\in U_+$ there exists a matrix    $D(i,j)\in Q $  such that $ (\tau(i))^{\top}D(i,j)\tau(j)>0 .$

 Consider a matrix
\begin{equation}\label{Tildex}\widetilde{D} := \displaystyle\sum_{(i,j) \in U_{+}}\alpha_{ij}D(i,j) + \displaystyle\sum_{j \in J_Q}\displaystyle\sum_{k \in
P\setminus M_Q(j)}\beta_{kj}A(k,j), \end{equation}
where the coefficients $\alpha_{ij}, (i,j) \in U_{+}$, and $ \beta_{kj}, k \in  P\setminus M_Q(j),  j\in J_Q, $ are such that
$$\displaystyle\sum_{(i,j) \in U_{+}}\!\alpha_{ij}+\displaystyle\sum_{j \in J_Q}\displaystyle\sum_{k \in  P\setminus M_Q(j)}\!\beta_{kj}=1,\;
  \alpha_{ij}>0 \; \forall (i,j) \in U_{+}, \ \beta_{kj} >0 \  \forall k \in  P\setminus M_Q(j),\; \forall j\in J_Q.$$
Notice that since the set $Q$ is convex, it holds $\widetilde{D} \in Q.$ Let us show that the following inequalities  are valid:
\begin{equation}\label{Des1}
 t^{\top}\widetilde{D} t>0 \ \forall t\in ({\rm conv\, Z_Q})\setminus T_0(Q);
\end{equation}

\vspace{-3mm}

\begin{equation}\label{Des2}
e^{\top}_k \widetilde{D} \tau(j)>0 \ \forall k\in P\setminus M_Q(j),\ \forall j\in J_Q.
\end{equation}
Inequalities (\ref{Des2}) are valid by construction.
 Moreover,  it holds (see (\ref{N*}))
 \begin{equation}\label{Cond}(\tau(i))^{\top}\widetilde{D} \tau(j) =0   \; \mbox{ if }    (i,j)\in U\setminus U_{+},  \;\;
 (\tau(i))^{\top} \widetilde{D} \tau(j)>0 \; \mbox{ if } (i,j)\in  U_{+}. \end{equation}
 To prove (\ref{Des1}), let us suppose that, on the contrary, there exists $\bar t\in ({\rm conv\, Z_Q})\setminus T_0(Q)$ such that

 \vspace{-5mm}

 \begin{equation}\label{EQ1}
 \bar t^{\top} \widetilde{D}\bar t=0.
\end{equation}
 From the condition $\bar t\in \rm conv\, Z_Q$, it follows that for some  set $\bar J\subset J_Q$, it holds
 %\begin{equation}\label{EQ2}
 $$\bar t = \displaystyle\sum_{j \in \bar J}{\alpha}_{j}\tau(j),\mbox{  where }  \displaystyle\sum_{j \in \bar J}{\alpha}_{j}=1, \ \alpha_{j}>0 \ \forall j\in \bar J.
$$%\end{equation}
   Then  we can rewrite (\ref{EQ1}) in the form
 $$0=\displaystyle\sum_{j \in \bar J}({\alpha}_{j}\tau(j))^{\top} \widetilde{D}\displaystyle\sum_{i \in {\bar J}}{\alpha}_{i}\tau(i)=\displaystyle\sum_{i \in {\bar J}}\displaystyle\sum_{j \in \bar J}{\alpha}_{i}{\alpha}_{j}(\tau(i))^{\top} \widetilde{D} \tau(j),$$
 wherefrom, taking into account (\ref{Cond}), we obtain
 $  (i,j) \in U\!\setminus \!U_{+} \, \forall i\in \bar J,$ $ \forall j\in \bar J,\;$ $ i\not=j. $ Therefore  $ (\tau(i))^{\top} D \tau(j)=0 \   \forall i\in \bar J, \forall j\in \bar J, \ \forall D\in Q.$
  Hence
 $$\bar t^{\top}D \bar t = \displaystyle\sum_{i \in \bar J}\displaystyle\sum_{j \in \bar J}{\alpha}_{i}{\alpha}_{j}(\tau(i))^{\top}D \tau(j)=0 \ \ \forall D\in Q,$$
 and  $\bar t\in {T_0}(Q)$, which  contradicts the condition $\bar t\in ({\rm conv\, Z_Q})\setminus {T_0}(Q)$.

 \vspace{2mm}

 Following the proof of   Lemma 3 {in} \cite{KT-SetValued}, it is easy to show  that  there exists $\bar A\in Q$, such  that
 \be t^{\top}\bar At >0 \ \forall t\in \Omega(Z_Q),\label{o-11}\ee
 where the set $ \Omega(Z_Q)$ is defined in (\ref{0**})-(\ref{set-hat1}) with $V=Z_Q$.

 Consider a matrix $$\bar D:=\frac12( \widetilde{D}+\bar A),$$  where $ \widetilde{D} \in Q$ is defined in (\ref{Tildex}) and $\bar A\in Q$ satisfies (\ref{o-11}).
   By construction, $ \bar D \in Q$ and  it holds
 \begin{equation}\begin{split}\label{8}&t^{\top}\bar D t >0\ \forall t\in \Omega(Z_Q)\cup ( ({\rm conv \,Z_Q})\setminus {T_0}(Q));\\
  &e^{\top}_k \bar {D} \tau(j)>0 \ \forall k\in P\setminus M_Q(j),\ \forall j\in J_Q.\end{split}\end{equation}

  To complete the proof of the theorem, let us show that {for the}  matrix  $\bar D$  defined above, the following inequalities   hold:
  \be  t^{\top}\bar D t>0\;\; \forall   t\in \{t\in T: \ \rho(t,{\rm conv} Z_Q)< \sigma(Z_Q)\}\setminus {\rm conv Z_Q}.\label{1234}\ee

  Suppose that, on the  contrary, there exists  a vector  $\bar t $, $$\bar t \in \{t\in T:  \rho(t,{\rm conv} Z_Q)< \sigma(Z_Q)\}\setminus {\rm conv Z_Q},$$ such that $  \bar{t}^{\ \top}\bar D\bar{t}=0.$ To come to a contradiction, we will take the next steps.

\vspace{3mm}

\textbf{Step {0}.} Set $s:= 1, \ t^s:=\bar t \; \in T.$

\textbf{Step 1.}  For a given  $t^s\in T $, if  either $t^s \in {\rm conv Z_Q}$ or $\rho(t^s, {\rm conv} Z_Q)\geq \sigma(Z_Q) $ or $(t^s)^{\top}\bar D t^s>0$, then STOP. Otherwise, i.e. if \begin{equation}\label{10} t^s\in \{t\in T: \rho(t^s, {\rm conv} Z_Q)< \sigma(Z_Q)\}\setminus {\rm conv\, Z_Q} \ \mbox{ and } (t^s)^{\top}\bar D t^s=0,
 \end{equation}  go to the next step.

 \textbf{Step 2.} As $\rho(t^s,{\rm conv }Z_Q)<\sigma(Z_Q),$ then, according to  Proposition \ref{Pnew0}, there exists an index $i_s\in J_Q$ such  that $P_0(t^s)\subset P_0(\tau(i_s))$, or equivalently, $P_+(\tau(i_s))\subset P_+(t^s)$.   Compute $ \theta_k=  \frac{t_k^s}{\tau_k(i_s)}>0$ $  \forall  k\in P_+(\tau(i_s))\subset P_+(t^s),$   and set $$\theta:=\min\{\theta_k, k\in P_+(\tau(i_s))\}>0.$$
 Let us show that $\theta <1$. If suppose that $\theta \geq 1,$ then, evidently,
\begin{equation}t_k^s\geq \tau_k(i_s) \ \forall k\in P_+(\tau(i_s)).\label{11-new}\end{equation} Taking into account the latter inequalities and the definition (\ref{SetT}) of the set $T$, we get
 \begin{equation}\label{11} 1=\displaystyle\sum_{k \in P_+(t^s)}t_k^s\geq \displaystyle\sum_{k \in P_+(\tau(i_s))}t_k^s\geq \displaystyle\sum_{k \in P_+(\tau(i_s))}\tau_k(i_s)= 1.
 \end{equation}
 From (\ref{11-new}) and (\ref{11}), it follows $t^s=\tau(i_s)$, which contradicts the assumption $t^s\notin {\rm conv \,Z_Q}.$ Therefore, we have shown that $0<\theta<1.$

 \vspace{2mm}

 Set $t^{s+1}:=(t^s-\theta \tau(i_s))/(1-\theta)$. Let us show that  the following relations are satisfied:
$$ \mbox{\bf a) }  t^{s+1} \in T; \ \ \mbox{\bf b) }  t^{s+1} \notin {\rm conv Z_Q};\ \ \mbox{\bf c) } (t^{s+1})^{\top} \bar  D t^{s+1}=0;\ \ \mbox{\bf d) } \rho(t^{s+1},{\rm conv} Z_Q)< \sigma(Z_Q).$$

The proof of these items is as follows.

\begin{enumerate}
\item[\bf{a)}]  By construction, $ t^{s}-\theta \tau(i_s) \geq 0.$ Hence $t^{s+1}\geq 0$ and $\mathbf{e}^{\top}t^{s+1}=1$. Therefore $t^{s+1}\in T.$
\item[\bf{b)}] Suppose that, on the contrary, $t^{s+1} \in {\rm conv Z_Q}.$ Then $t^s=(1-\theta)t^{s+1}+\theta \tau(i_s), $ where $\theta \in (0, 1), $   $\tau(i_s) \in {\rm conv Z_Q},$  and $t^{s+1}\in {\rm conv Z_Q}.$  Hence $t^s \in {\rm conv Z_Q}$ that contradicts (\ref{10}). The obtained contradiction  proves that $t^{s+1} \notin {\rm conv Z_Q}.$
\item[\bf{c)}] Since $t^{s+1}\in T,$ then $(t^{s+1})^{\top} \bar D t^{s+1}\geq 0$ and the following relations hold:
\begin{equation}\begin{split}0\leq (t^{s+1})^{\top}\bar D t^{s+1}&=(1-\theta)^{-2}(t^s-\theta (\tau(i_s))^{\top} \bar D(t^s-\theta (\tau(i_s)))\\
&=(1-\theta)^{-2}(-2(t^s)^{\top}\bar D\tau(i_s))\leq 0.\label{10-new}\end{split}\end{equation}
Here we took  into account that $t^\top D \tau(j)\geq 0\;\; \forall  t \in T, \forall j \in J_Q,\forall D\in Q.$

From (\ref{10-new}), it follows: $(t^{s+1})^{\top} \bar D t^{s+1}= 0.$
\item[\bf{d)}] Suppose that $\rho(t^{s+1},{\rm conv} Z_Q)\geq \sigma(Z_Q).$
From the condition {\bf a)} and  the inequalities  (\ref{8}), we get
 $(t^{s+1})^{\top} \bar D t^{s+1}>0$,  which contradicts {\bf c)}.
  Hence $\rho(t^{s+1},{\rm conv} Z_Q)< \sigma(Z_Q)$.
\end{enumerate}

    Notice that, by construction, the number of the null components of the  vector   $t^{s+1}$ is larger than that of the  vector  $t^{s}$:
    \begin{equation}\label{12} |P_0(t^{s+1})|\geq |P_0(t^{s})|+1.\end{equation}

Let us substitute $s$ by $s+1$ and go to  \textbf{Step 1}.

\vspace{3mm}

  On the \textbf{Step 1}, the situation STOP cannot  happen   since  for any $s\geq 0$, by construction,  relations (\ref{10}) hold.
 On another hand, due to (\ref{12}), one cannot repeat the described above procedure more than $p$ times.
   This contradiction proves that inequalities (\ref{1234}) hold true. The theorem is proved. $  \blacksquare$

\vspace{4mm}

 It follows from relations (\ref{For2}) that the constructed in the proof of Theorem \ref{L25-2-2} matrix $ \bar D$ is a minimally active element of the set $Q$ in the case ${T_0}(Q)\not =\emptyset$.
  If ${T_0}(Q)=\emptyset$, then the set $Q$ satisfies the Slater condition, and, hence, there exists a matrix  $\widehat D \in Q $ such that ${T_0}(\widehat D)=\emptyset. $ The matrix $\widehat D$ is a minimally active element in   this case.

\begin{corollary}\label{Lev} Under the conditions of Theorem \ref{L25-2-2}, there exists $\bar D\in Q$ such that
\be  \ t^\top  \bar D t>0\; \ \forall t \in \Omega(Z_Q); \ e^\top_k \bar D\tau(j)> 0\ \ \forall k \in P\setminus  M_Q(j),\ \forall  j \in J_Q.\label{2in11}\ee
\end{corollary}

\section{On equivalent representations of a  face    of the cone ${\cop}$ }\label{S5}

In this and the subsequent sections, based on the results above, we will obtain  equivalent representations of   faces of the cone ${\cal COP}^p$ and  their dual ones which will be useful for  further research and  some applications (for example, for regularization procedures  and a study of the facial structure of $\cop$).

\vspace{2mm}

Given  a finite non-empty   vector set $V$ in the form (\ref{0**}), $V=\{ t(i) \in T , \ i\in I \}$, and a set ${\cal L}:=\{L(i), \ i\in I\}$  of sets
$L(i)$  such that  $ P_+(t(i)) \subset L(i)\subset P \ \forall i \in I,$
consider a cone
\be {\mathcal K}={\mathcal K}(V, {\cal L}):
=\{D\in {\cal COP}^{p}:\ e^\top_kDt(i)=0 \ \forall k \in L(i), \ \forall   i \in I\}.\label{cone01}\ee

It follows from the definition of the cone $\mathcal K$ that
$(t(i))^\top Dt(i)=0$  for all $i \in I$ and $D\in {\cal K}$. Consequently, for all $D\in {\cal K}$, the vectors  $t(i) , i\in I,$ are optimal solutions of the problem (\ref{1problem}).
%\be \min \; t^\top Dt\;\; {\rm s.t. } \; t \in T.\label{problem}\ee
Hence
\be e^\top_kDt(i)=0 \ \ \forall k \in L(i); \ \;  e^\top_kDt(i)\geq 0 \ \ \forall  k \in P\setminus L(i), \;\forall  i \in I, \; \forall D \in {\cal K},\label{in}\ee
and the cone ${\cal K}$ defined in (\ref{cone01}) can be rewritten in the form
\begin{equation}\begin{split}  {\mathcal K}=\widetilde{ \ck}:
=\{D\in &{\cal S}(p): \  t^\top Dt\geq 0 \; \forall t \in T;\\&e^\top_kDt(i)=0 \ \forall k \in L(i), \; e^\top_kDt(i)\geq 0\ \ \forall k \in P\setminus L(i), \ \forall i \in I\},\label{cone1}\end{split}\end{equation}
\begin{proposition}\label{P-00}
The  cone ${\mathcal K}$ defined in (\ref{cone01}) is a face of ${\cal COP}^p.$
\end{proposition}
{\bf Proof.} Let
\be A\in {\cal COP}^p, \; B\in {\cal COP}^p  \;\mbox{ such that} \ (A+B)\in {\cal K}.\label{ZZ}\ee
 Remind that by definition, ${\cal K}$ is a face of ${\cal COP}^p$ if  relations (\ref{ZZ})  imply $A \in {\cal K}$ and $B\in {\cal K}.$
By {construction}, the condition $( A+B)\in {\cal K}$  is equivalent to the conditions
\be ( A+B)\in {\cal COP}^p, \; e^\top_k( A+B)t(i)=0\ \ \forall k \in L(i),\ \forall i \in I,\label{ZZ1}\ee
which imply the equalities
 $$(t(i))^\top ( A+B)t(i)=0, \;  (t(i))^\top  At(i)=0, \; ( t(i))^\top Bt(i)=0 \ \forall i \in I.$$
Moreover,  from the conditions
$ A\in {\cal COP}^p,$ $ (t(i))^\top At(i)=0$, and $B\in {\cal COP}^p,$ $ (t(i))^\top Bt(i)=0$,  it  follows  that
$  At(i)\geq 0$ and $ Bt(i)\geq 0\  \ \forall  i\in I.$ Taking into account these inequalities  and the equalities in (\ref{ZZ1}), we obtain $$
e^\top _kAt(i)=0,\; e^\top _kBt(i)=0 \ \forall  k \in L(i),\  \forall  i \in I.$$
By the definition of the cone ${\cal K}$, from the last equalities and relations (\ref{ZZ}), it follows that $A\in {\cal K}$ and $B\in {\cal K}$. The proposition is proved. $ \blacksquare$
\begin{remark}
Note that if in (\ref{cone01}) we have $L(i)=P_+(t(i))\  \ \forall  i \in I, $ then the cone $\ck$ is an exposed face of $\cop$  (for the definition of the exposed face see e.g. \cite{Pataki}).
\end{remark}\begin{remark}
 In  what follows (see Corollary \ref{Col-2} in section \ref{S7}), we will show that any face of the cone $\cop$ can be presented in the form (\ref{cone01}).
\end{remark}

Given  the cone ${\cal K}$  defined in (\ref{cone01}), consider
the set $  {T_0}(\ck)$ of all its {zeros},
the  set $ Z_\ck=\{\tau(j), j \in J_\ck\} $ of  all vertices of the set ${\rm conv} T_0(\ck)$  (i.e. the set of minimal zeros of $\ck$),
and  the sets $  M_\ck(j), j \in J_\ck,$ defined in (\ref{Formula2}), (\ref{vert}), and (\ref{M}) with $Q=\ck.$
Evidently, for all $i\in I,$  it holds
 $t(i)\in {T_0}({\cal K}).$

Note that in general, there may exist a vector $\bar {t} \in {T_0}({\cal K})$ such that $\bar {t} \not \in {\rm conv} \{t(i), i \in I\}.$  In fact, suppose that $V=\{t(1), \; t(2)\}$ and $P_+(t(1))\subset P_+(t(2)).$ For a sufficiently small $\theta>0$, consider
$t(\theta):=(t(2)-\theta t(1))/(1-\theta).$
By construction, $t(\theta)\in T$ and $t(\theta)\not\in {\rm conv} \{t(1), t(2)\}.$
From the definitions of the sets $P_+(t),$ $P_0(t),$  and  the cone ${\cal K}$, and from the assumption $P_+(t(1))\subset P_+(t(2))$,  it  follows :
$$(t(i))^\top Dt(i)=0,\; i=1,2; \;\;( t(1))^\top Dt(2)=0 \;\forall  D\in {\cal K}.$$
Then $(t(\theta))^\top Dt(\theta)=0 \;\forall D\in {\cal K}$ and, hence, $t(\theta) \in T_0({\cal K}).$

\begin{lemma} \label{OL-3}  The defined in  (\ref{cone01}) cone ${\mathcal K}$ coincides with the following cone $\bar{\mathcal K}$:
\be \begin{split} \bar{\mathcal K}:=\{D\in & {\cal S}(p): \
t^\top Dt\geq 0 \; \forall t \in \Omega(Z_\ck); \\
& e^\top_kD\tau(j)=0 \  \forall  k \in M_{\ck}(j), \; e^\top_kD\tau(j)\geq 0\  \forall  k \in P\setminus  M_{\ck}(j), \   \forall  j \in J_\ck\},\end{split}\label{bar-cone}\ee
and there exists a matrix $\bar D\in \bar{\mathcal K} $ such that
\be  \ t^\top\bar Dt>0\; \ \forall t \in \Omega(Z_\ck), \; e^\top_k \bar D \tau(j)> 0\  \forall  k \in P\setminus  M_{\ck}(j),  \  \forall  j \in J_\ck,\label{in11}\ee
where the sets  $Z_\ck=\{\tau(j),\in J_\ck\},$  $M_\ck(j), j \in J_\ck,$ and  $\Omega(Z_\ck)$ are defined in (\ref{vert}), (\ref{M}) with $Q=\ck$ and in (\ref{o-1}),(\ref{set-hat1}) with $V=Z_\ck$.
\end{lemma}
{\bf Proof.}
 First, let us show that the defined in (\ref{cone01})  cone ${\cal K}$  coincides with the cone
\begin{equation}\begin{split}  \widehat{\cal K}:=\{ D\in &\mathcal S(p): \  t^\top Dt\geq 0\;\forall t \in T; \\
&e^\top _kD\tau(j)=0 \  \forall  k\in M_{\ck}(j),\;e^\top_k D\tau(j)\geq 0\  \forall  k\in P\setminus M_{\ck}(j),\   \forall    j \in J_\ck\}.\label{in12}\end{split}\end{equation}
 From the definitions  of the  sets $\ck$ and ${T_0}(\ck)$, it follows:
$$\ck\subset \,\cop  \mbox{ and }  \; (\tau(j))^\top D\tau(j)=0 \  \forall  j\in J_\ck,\; \forall D\in \ck.$$
Then, for all $D\in \ck$,  the vectors $\tau(j), j \in J_{\ck},$ are  optimal solutions of problem (\ref{1problem}). Hence relations (\ref{1in}) with $Q=\ck$ hold true.
Moreover, from  (\ref{M}) with $Q=\ck$  it follows:
$$e^\top_kD\tau(j)=0 \  \ \forall  k \in M_{\ck}(j),  \ \forall   j \in J_\ck, \;{\forall} D\in \ck.$$
These relations together with (\ref{1in}) imply the inclusions $\ck \subset \widehat{\ck}$
  and $P_+(\tau(j))\subset M_{\ck}(j)$ $\  \forall  j \in J_\ck.$

Now,  consider any vector $t(i)\in V$. Since, by construction, $t(i)\in {T_0}(\ck)$, then there exists a subset $J(i)\subset J_\ck$ and numbers $\alpha_j=\alpha_j(i), j \in J(i),$ such that
\be t(i)=\sum\limits_{j\in J(i)}\alpha_j\tau(j),\; \sum\limits_{j\in J(i)}\alpha_j=1,\; \alpha_j>0  \ \forall  j \in J(i).\label{tau}\ee
Hence, by construction,
$$0=e^\top_kDt(i)=\sum\limits_{j\in J(i)}\alpha_je^\top_kD\tau(j) \ {\forall} k\in L(i), \; \forall D\in \ck.$$
Taking into account the equalities above and  relations (\ref{1in}) (with $Q=\ck$), we conclude that
$e^\top_kD\tau(j) =0\  \forall j \in J(i), \;\forall k\in L(i), \; \forall D\in \ck$, wherefrom it comes
\be L(i)\subset M_{\ck}(j)\ \forall j \in J(i), \; \; \forall i \in I.\label{LM}\ee

Now, we will show that $\widehat{\ck}\subset \ck.$ Let $\bar D\in \widehat{\ck}$ and $t(i) \in V.$ It follows from (\ref{tau}) that
$$e^\top_k\bar Dt(i)=\sum\limits_{j\in J(i)}\alpha_je^\top_k\bar D\tau(j)\ \forall k \in P.$$
Taking into account these equalities and inclusions (\ref{LM}), we conclude that
$e^\top_k\bar Dt(i)=0 \;\forall  k \in L(i),$ $\forall i \in I$, and hence, $\bar D\in \ck$. Thus we have shown that $\widehat{\ck}\subset \ck.$  The equality
$\ck=\widehat{\ck}$  is proved.

From the definition of $\widehat{\ck}$  (see (\ref{in12})), it follows  $\widehat{\cal K} \subset \bar\ck$ and from Theorem \ref{LO-1}  we get $\bar\ck \subset \widehat{\cal K}. $ These inclusions together with  the equality $\ck=\widehat{\ck}$ imply  $\ck=\bar\ck.$

\vspace{1mm}

Applying   Corollary \ref{Lev}   with $Q=\ck$, we conclude that there exists   a matrix  $\bar D\in \bar{\mathcal K}=\ck $ such that relations
(\ref{in11}) hold true. The lemma is proved. $  \blacksquare$

\vspace{2mm}

It is worth noting that for  the cone ${\cal K}$ defined in (\ref{cone01}), the obtained in this section  equivalent  representations (\ref{cone1}), (\ref{bar-cone}), and  (\ref{in12})  may    be  more  preferable than  the original definition, particularly (see the next section) when one needs to describe its dual cone.

\section{Alternative representations of the dual cone to the  defined in (\ref{cone01}) face of $\cop$}\label{S6}

In this section, we will use  the following statements proved in \cite{Rock} (See Theorem 6.5 and Proposition 16.4.2)
\begin{proposition} \label{Ro} For two closed convex cones $C_1$ and $C_2$ in $\mathbb R^m$,  it holds true
\be (C_1 \cap C_2)^*={\rm  cl }(C_1^*
 \oplus  C_2^*).\label{RR}\ee

 \vspace{-2mm}

 If ${\rm relint} (C_1) \cap {\rm relint}  (C_2) \not =\emptyset$ , then

\vspace{2mm}

(i) $(C_1^*  \oplus  C_2^*)$ is a closed set
and the closure operation in (\ref{RR}) can be omitted;

(ii) ${\rm relint}(C_1\cap C_2) = {\rm relint}(C_1) \cap  {\rm relint}(C_2).$
\end{proposition}
 Here and in what follows, $C^*$ denotes the dual cone for a cone $C$ and $\oplus $ denotes the  Minkowski  sum.

Based on this proposition and   definition (\ref{cone01}), the  dual cone   to the cone ${\cal K} $   can be written in the form
\be \label{dual01} {\cal K}^*={\rm cl}\, \cal G,\ee
where

 \vspace{-10mm}

\begin{equation*}\begin{split}\mathcal G:=&\{D\in {\cal S}(p): \ D=\sum\limits_{i=1}^{p_*}\alpha_i\mu(i)(\mu(i))^\top+\sum\limits_{i\in I}(\lambda(i)(t(i))^\top+t(i)(\lambda(i))^\top),\\ & \mu(i)\in T, \alpha_i\geq 0,i=1,...,p_*;\;  \lambda_k(i)=0\  \forall  k \in P\setminus L(i),\  \forall  i \in I\},\end{split}\end{equation*}
$ p_*=p(p+1)/2.$  Let us give alternative descriptions of the dual cone $\ck^*.$

It follows from Lemma \ref{OL-3} that the cone ${\cal K}$ admits the following representation:
\be {\cal K}={\cal COP}(Z_\ck)\cap {\cal K}_{pol},\label{pol}\ee
where
\be {\cal COP}(Z_\ck):=\{D\in {\cal S}(p):\ t^\top Dt\geq 0 \;\; \forall t \in \Omega(Z_\ck)\},\qquad\qquad\label{12-4}\ee
\begin{equation*}\begin{split}{\cal K}_{pol}:=\{D\in {\cal S}(p): \ &e^\top_kD\tau(j)= 0\  \ \forall  k \in  M_{\ck}(j), \\ &e^\top_kD\tau(j)\geq 0 \ \forall   k \in P\setminus  M_{\ck}(j) , \  \forall j \in J_\ck\},\end{split}\end{equation*}
 the sets  $Z_\ck=\{\tau(j),\in J_\ck\},$  $M_\ck(j), j \in J_\ck,$ and  $\Omega(Z_\ck)$ are defined in (\ref{vert}), (\ref{M}) with  $Q=\ck$ and
(\ref{o-1}), (\ref{set-hat1}) with $V=Z_\ck$.

It is known (see \cite{E-P})   that the cone ${\cal COP}(Z_\ck)$ is convex, closed, and pointed, its interior can be presented in the form
$${\rm int}({\cal COP}(Z_\ck))=\{D\in {\cal S}(p):t^\top Dt>0 \; \forall t \in  \Omega(Z_\ck)\}$$
and  its dual cone  is $({\cal COP}(Z_\ck))^*={\rm cl}\, {\cal G}(  Z_\ck),$ where
$${\cal G}(Z_\ck):=\left \{D\in {\cal S}(p){: }\ D= \sum\limits_{i=1}^{p_*}\alpha_i\mu(i)(\mu(i))^\top,\ \alpha_i\geq 0,\; \mu(i)\in \Omega(Z_\ck) \ \forall  i=1,...,p_*\right\}.$$

By construction, the set $\Omega(Z_\ck) $ is closed and   ${\mathbf e}^\top \mu=1,$ $\mu \geq 0$ for any $\mu \in \Omega(Z_\ck)$. Then we can show that $({\cal COP}(Z_\ck))^*={\cal G}(Z_\ck).$

 The cone ${\cal K}_{pol}$ is convex, closed, and for this cone it holds
$${\rm relint}\,{\cal K}_{pol}\!\supset\! \{D\in {\cal S}(p)\!:\!e^\top_kD\tau(j)= 0\, \forall  k \in M_{\ck}(j);  \ e^\top_kD\tau(j)> 0\, \forall  k \in P\setminus  M_{\ck}(j),\  \forall j \in J_\ck\},$$
\begin{equation*}\begin{split}({\cal K}_{pol})^*=\{D\in {\cal S}(p):   D=\sum\limits_{j\in J_\ck}&(\lambda(j)(\tau(j))^\top+
\tau(j)(\lambda(j))^\top), \\ & \lambda_k(j)\geq 0\  \forall  k \in P\setminus M_{\ck}(j),\  \forall  j\in J_\ck\}.\end{split}\end{equation*}
Then it follows from {Lemma} \ref{OL-3} that there exists a matrix $\bar D$ such that
$$\bar D\in {\rm int}({\cal COP}(Z_\ck))\cap {\rm relint}\,{\cal K}_{pol},$$%\label{int}\ee
and based on Proposition \ref{Ro} we conclude that
\begin{equation*}\begin{split} {\rm relint}({\cal K})&={\rm relint}({\cal COP}(Z_\ck)) \cap  {\cal K}_{pol}) = {\rm int}({\cal COP}(Z_\ck))) \cap  {\rm relint}({\cal K}_{pol}),\\
\ck^*&=({\cal COP}(Z_\ck) \cap  {\cal K}_{pol})^* = ({\cal COP}(Z_\ck))^*\oplus ({\cal K}_{pol})^*.\end{split}\end{equation*}

  Thus  we have proved the following  theorem.

 \begin{theorem}
  Let the cone ${\cal K}$ be  defined in (\ref{cone01}). Then its dual cone can be described as   ${\cal K}^*=\bar{{\cal G}}$, where
 \begin{equation}\begin{split}\bar {\cal G}:=&\bigl\{{D\in} {\cal S}(p): D= \sum\limits_{i=1}^{p_*}\alpha_i \mu(i)(\mu(i))^{\top}+\sum\limits_{j\in J_\ck}(\lambda(j)(\tau(j))^\top+\tau(j)(\lambda(j))^\top),\label{dualcone}\\ &  \alpha_i\geq 0,\; \mu(i)\in \Omega(Z_\ck)\  \forall  i=1,...,p_*;\ \lambda_k(j)\geq 0 \  \forall  k \in P\setminus M_{\ck}(j),\  \forall  j \in J_\ck\bigl\},\end{split}\end{equation}
 and the sets  $Z_\ck=\{\tau(j),\in J_\ck\},$  $M_\ck(j), j \in J_\ck,$ and  $\Omega(Z_\ck)$ are defined in (\ref{vert}), (\ref{M}) with $Q=\ck$, and in (\ref{o-1}), (\ref{set-hat1}) with $V=Z_{\ck}$.
 \end{theorem}

Notice that in the case when the information about the  minimal zeros  $\tau(j)$ and sets $M_{\ck}(j)$, $j \in J_\ck,$  is not available, one can use  an alternative  representation of the dual cone ${\cal K}^*$,
\be{\cal K}^*={\rm cl}\,\widetilde{\cal G}, \label{dualcone2}\ee
where
\begin{equation}\begin{split} \widetilde{\cal G}:=\bigl\{D\in & {\cal S}(p):\ D=\sum\limits_{i=1}^{p_*}\alpha_i\mu(i)(\mu(i))^{\top}+\sum\limits_{i\in I}(\lambda(i)(t(i))^\top+t(i)(\lambda(i))^\top),\\ & \alpha_i\geq 0,\; \mu(i)\in T\  \forall   i=1,...,p_*;\; \lambda_k(i)\geq 0 \  \forall  k \in P\setminus L(i), \  \forall  i \in I \bigl\}.\end{split}\end{equation}
This description of the dual cone ${\cal K}^*$ is based on its  representation (\ref{cone1})  and formula (\ref{RR}).

The obtained here representations (\ref{dualcone}) and (\ref{dualcone2}) of the dual cone ${\cal K}^*$ are  preferable than formula (\ref{dual01}) directly based  on  the definition (\ref{cone01}) of the cone ${\cal K} $. This can be motivated  by the following reasons:
\begin{itemize}
\item the inclusions
${\cal G}\subset \widetilde{\cal G}\subset \bar {\cal G}$ hold true
and, in general,  ${\cal G}\not=\widetilde{\cal G}$;
\item  to construct the set $\widetilde{\cal G}$, we   need only the original  data   used in the definition (\ref{cone01}) of the cone $\ck$;
\item in the representation ${\cal K}^*=\bar{\cal G}$,  the closure operator  is absent.
\end{itemize}

In fact,  it is evident that ${\cal G}\subset \widetilde{\cal G}$.

Let us show that $\widetilde{\cal G}\subset \bar {\cal G}.$ Suppose that $D \in \widetilde{\cal G}$, wherefrom by definition,
\be D=\sum\limits_{i\in I_*}\alpha_i\mu(i)(\mu(i))^{\top}+\sum\limits_{i\in I}(\lambda(i)(t(i))^{\top}+t(i)(\lambda(i))^\top),\label{DZ}\ee
 where  $\mu(i)\in T,\alpha_i>0 $  $ \ \forall i\in I_*\subset \{1,...,p_*\}$,  $  \lambda_k(i)\geq 0 \ \forall  k \in P\setminus L(i), \forall i \in I.$

Since, by construction $t(i)\in {\rm conv}\, Z_\ck,\ i \in I,$  there exist numbers $\alpha_{ij}\geq  0,\, j\in J_\ck$, $i\in I,$ such that
$t(i)=\sum\limits_{j \in J_\ck}\alpha_{ij}\tau(j)\   \forall  i \in I.$
Taking into account these equalities, we can rewrite equality (\ref{DZ}) as follows:
\be D=\sum\limits_{i\in I_*}\alpha_i\mu(i)(\mu(i))^{\top}+\sum\limits_{j\in J_\ck}(\bar \lambda(j)(\tau(j))^{\top}+\tau(j)(\bar\lambda(j))^\top),\label{2DZ}\ee
where $\bar \lambda(j):=\sum\limits_{i\in I}\alpha_{ij}\lambda(i),$  $j \in J_\ck.$
Let us show that
\be \bar \lambda_k(j)\geq 0\  \ \forall  k \in P\setminus M_\ck(j), \ \forall  j \in J_\ck.\label{10-1}\ee
Suppose that on the contrary, there exist $j_0\in J_\ck$ and $k_0\in P\setminus M_\ck(j_0)$, such that
$\bar \lambda_{k_0}(j_0)=\sum\limits_{i\in I}\alpha_{ij_0}\lambda_{k_0}(i)<0.$ Hence, there exists $i_0\in I$ such that
$\alpha_{i_0j_0}>0,$ $\lambda_{k_0}(i_0)<0.$ This implies $j_0\in J(i_0):=\{j \in J_\ck:\alpha_{i_0j}>0\},$ $ k_0\in L(i_0).$ From these relations and (\ref{LM}), it follows that $k_0\in M_\ck(j_0).$ But this contradicts the condition  $k_0\in P\setminus M_\ck(j_0)$. Inequalities (\ref{10-1}) are proved.

 If $\mu(i)\in \Omega(Z_\ck)$ for all $i \in I_*$, then  from (\ref{2DZ}) and (\ref{10-1}), it  follows $D\in \bar {\cal G}.$
 Suppose, first, that there exists $i_0\in I_*$ such that $0<\rho(\mu(i_0), {\rm conv}\, Z_\ck)<\sigma(Z_\ck).$  Replacing  $\mu(i_0)$ by $\bar {\mu}(i_0):=\beta \mu(i_0)$ and
 $\alpha_{i_0}$ by $\bar \alpha_{i_0}:=\alpha_{i_0}/\beta^2$ with  $\beta:=\sigma(Z_\ck)/\rho(\mu(i_0), {\rm conv}Z_\ck)>0$, we obtain
 \be \alpha_{i_0}\mu(i)(\mu(i_0))^\top=\bar\alpha_{i_0}\bar\mu(i_0)(\bar\mu(i_0))^\top,  \;\bar {\mu}(i_0)\in \Omega(Z_\ck).\label{10-2}\ee
  Now suppose that there exists $i_0\in I_*$ such that $\rho(\mu(i_0), {\rm conv} Z_\ck)=0$. Hence
  $$\mu(i_0)=\sum\limits_{j \in \bar J}\beta_j\tau(j), \; \beta_j>0\  \ \forall  j \in \bar J\subset J_\ck.$$
  From these relations we conclude that
  \begin{equation}\begin{split}\label{10-3} &P_0(\mu(i_0))\subset P_0(\tau(j))\  \ \forall  j \in \bar J, \\
   & \alpha_{i_0}\mu(i_0)(\mu(i_0))^\top=\frac{\alpha_{i_0}}{2}\sum\limits_{j \in \bar J}\beta_j(\tau(j)({\mu}(i_0))^\top+ {\mu}(i_0)(\tau(j))^\top).\end{split}\end{equation}
  Set
  $\bar \alpha_{i_0}=0, $ $ \tilde  \lambda(j)=\bar \lambda(j)+\alpha_{i_0}\beta_j{\mu}(i_0), j \in \bar J, $   and  $ \tilde  \lambda(j)={\bar \lambda(j)}, i \in J_\ck\setminus \bar J.$ \\
  From (\ref{2DZ})-(\ref{10-3}), it follows  that the matrix  $D \in \widetilde{\cal G}$ can be written  in the form (\ref{dualcone}) with \\ $\qquad \alpha_i, \ \mu(i), i \in I_*,$ and $\lambda(j), j \in J_\ck, \ $ replaced by $\ \ \bar \alpha_i,\ \bar{\mu}(i) , i \in I_*;$    $ \tilde  \lambda(j), j \in J_\ck,\ $  such that   $\ \tilde \lambda_k(j)\geq 0 $ $  \ \forall  k \in P\setminus M_\ck(j), $ $ \forall j \in J_\ck,$  and
    $\  \bar \alpha_i>0$  $ \Rightarrow  {\bar \mu(i)}\in \Omega(Z_\ck)$  $\forall i\in I_*$. Therefore, we conclude that $ D\in \bar {\cal G}$ and, hence,  $\widetilde{\cal G}\subset \bar {\cal G}.$

   \vspace{4mm}

   To show that, in general, ${\cal G}\not= \widetilde{\cal G}$, let us consider an example.

    Let $p=2,$ $I=\{1\}$, $t(1)=(1,\; 0)^\top,$ $L(i)=\{1\},$ $\lambda(1)=(0,\; 1)^\top. $  Then
   { \small $$D_*=\lambda(1)(t(1))^\top+t(1)(\lambda(1))^\top=\left(\begin{array}{cc}0 &1\cr 1&0\end{array}\right)\in \widetilde{\cal G}.$$}
    It is easy to see that the matrix $D_*$ cannot be presented in the form
    \bea D_*=\sum\limits_{i=1}^3\alpha_i\mu(i)(\mu(i))^{\top}+
    \bar \lambda(1)(t(1))^{\top}+t(1) (\bar\lambda(1))^\top,\nonumber\eea
    where $\mu(i)\in T, \alpha_i\geq 0, i=1,2,3;$ $\bar \lambda(1)=(\bar \lambda_1(1),0).$ Hence $D_*\not \in {\cal G}.$

\section{The minimal face of ${\cal COP}^p$ containing a given convex set}\label{S7}

Let $Q$ be a convex closed subset of ${\cal COP}^p$ with the corresponding sets   $T_0(Q)$ and $Z_Q=\{\tau(j), j \in J_Q\} $ of all  zeros  and minimal zeros  of $Q$ defined in (\ref{Formula2}) and (\ref{vert}),  and  the sets $ M_Q(j), j \in J_Q,$ defined by relations (\ref{M}).

Consider the cone
$$ {\cal K}_Q:=\{D\in {\cal COP}^{p}:{e^\top_k}D\tau(j)=0 \ \forall  k \in M_Q(j), \  \forall   j \in J_Q\}.$$%\label{K-Q}\ee
It is evident that $Q\subset  {\cal K}_Q$ and it is easy to show that
$$ T_0 (\ck_Q)={T_0}(Q),\; Z_{\ck_Q}=\ Z_Q,\; M_{\ck_Q}(j)=M_Q(j) \;  \forall  j \in J_Q=J_{\ck_Q}.$$
 Then  from the results of sections \ref{S4} and \ref{S5}, one can conclude that ${\cal K}_Q$ is a face of ${\cal COP}^{p}$ and
 \be {\rm relint}({\cal K}_Q)={\rm relint}({\cal COP}(Z_Q) \cap  {\cal K}^Q_{pol}) = {\rm relint}({\cal COP}(Z_Q)) \cap  {\rm relint}({\cal K}^Q_{pol}),\label{relKA}\ee
 where  the sets  $\Omega(Z_Q)$ and   ${\cal COP}(Z_Q)$ are
 defined by the rules (\ref{o-1}),   (\ref{set-hat1}),
 and (\ref{12-4})
 using the  minimal zeros   set $Z_Q$ and
 $${\cal K}^Q_{pol}:=\{D\in {\cal S}(p):{e^\top_k}D\tau(j)=0 \  \forall  k \in M_Q(j);  \  {e^\top_k}D\tau(j)\geq 0  \  \forall  k \in P\setminus M_Q(j), j \in J_Q\}.$$

It follows from (\ref{2in11})  and (\ref{relKA}) that there exists a matrix  $\bar D$  such that
 \be  \bar  D \in {\rm relint}({\cal K}_Q)\cap Q\not=\emptyset.\label{QQW}\ee

The  following   Proposition is proved in \cite{Pataki2} (see  Proposition 3.2.2).

\begin{proposition}\label{Pna} {L}et $ F $ be a  face of a convex cone $ K$ and $Q $ a
convex subset of $K$.  If $Q\subset F$ and
 $Q\cap {\rm relint} F\not=\emptyset $ then $F = {\rm face}(Q,K)$.
\end{proposition}

Here  and in what follows ${\rm face}(S,C)$ denotes  the minimal (by inclusion) face of a cone $C$ containing a set $S.$

\begin{theorem}\label{T-face}
Let $Q$ be a convex closed subset of ${\cal COP}^p$. Then
$${\rm face}(Q,\cop)={\cal K}_Q:=\{D\in {\cal COP}^p:{e^\top_k}D\tau(j)=0 \ \forall  k \in M_Q(j), \  \forall  j \in J_Q\},$$
where
$\{\tau(j), j \in J_Q\} $ is the set of  all  minimal zeros of $Q$   and the sets
$ M_Q(j),  j \in J_Q,$ are defined in (\ref{M}).
\end{theorem}
{\bf Proof.} The statement of this theorem follows from condition
(\ref{QQW}) and Proposition \ref{Pna}.

\begin{corollary}\label{Col-2} Any face of $\cop$ can be presented in the form (\ref{cone01}) with some  vectors $t(i)\in T$ and sets $L(i)$, $P_+(t(i))\subset L(i)\subset P,$  $\  \forall  i \in I,$ $0\leq |I|<\infty.$
\end{corollary}

{\bf Proof.} Let $F$ be a face of the cone ${\cal COP}^p$. It is known that $F$ is a convex closed subset of ${\cal COP}^p$.  Applying Theorem
\ref{T-face} with $Q=F$, we obtain
$${\rm face}(F,\cop)=\{D\in {\cal COP}^p:{e^\top_k}D\tau(j)=0 \  \forall  k \in M_F(j),  \forall  j \in J_F\},$$
where $\{\tau(j), j \in J_F\} $ is
the set of  all   minimal zeros of $F$  and the sets
$ M_F(j),  j \in J_F,$ are defined in (\ref{M}) with $Q=F$.
Taking into account the evident equality ${\rm face}(F,\cop)=F$, we conclude that $F$ can be represented in the form (\ref{cone01}), namely
$$F=\{D\in {\cal COP}^p:{e^\top_k}D\tau(j)=0 \  \forall  k \in M_F(j), \  \forall  j \in J_F\}.$$ The corollary is proved. $  \blacksquare$

\vspace{2mm}

Consider a matrix $A\in \cop $   and   the corresponding set   of all its   minimal zeros $ {Z}_A=\{\tau(j), j \in J_A\}$.

{\begin{corollary} The minimal  face of $\cop$  containing a given copositive matrix $A$ is as follows:
$${\rm face}(A,\cop)=\{D\in \cop: e^\top_kD\tau(j)=0\;  \forall   k \in M_A(j), \  \forall  j \in J_A\},$$
where $M_A(j)=\{k \in P: e^\top_kA\tau(j)=0\},$ $ j \in J_A$.
\end{corollary}}

Notice that this corollary correlates with   results from \cite{Dick}, where for a given matrix $A\in {\cal COP}^p$, the author gives explicit description  for
${\rm span} \{{\rm face}(A,\cop)\} $.

\section{On equivalent descriptions of the  feasible set of a copositive problem}\label{S3}

Consider a convex copositive problem in the form
\begin{equation} \displaystyle\min_{x \in \mathbb R^n } \ c(x) \;\;\;
\mbox{s.t. }   {\mathcal A}(x)
\in \mathcal{COP}^p,\label{cop}\end{equation}
 where $x=(x_1,...,x_n)^{\top}$  is the vector of  decision variables,  $c: \; \mathbb R^n\to \mathbb R$ is a convex function, and $\mathcal A: \; \mathbb R^n\to {\cal S}(p)$ is a given matrix  function
 such that for any  $x\in \mathbb R^n, y \in \mathbb R^n, $ the following inclusions are satisfied:
\begin{equation}\label{A}{\cal A}(\lambda x+(1-\lambda)y)-\lambda {\cal A}(x)-(1-\lambda){\cal A}(y)  \in \cop\;\; \forall \, \lambda\in [0,1].\end{equation}

The   aim of this section is to  prove Theorem \ref{L25-02-1} (see below) and  use it to obtain equivalent useful  descriptions of the feasible set of problem   (\ref{cop}). {This set can be written in the form
\be X:=\{x\in \mathbb R^n: {\mathcal A}(x)\in {\cal COP}^p\}
=\{x\in \mathbb R^n :t^{\top}{\mathcal A}(x)t\geq 0 \;\; \forall  t \in T\},\label{setX}\ee
where  $\mathcal A(x)$ is a matrix function  satisfying (\ref{A}) and the set $T$ is defined in (\ref{SetT}).

Let $\cd$ be a subset of the cone $\cop$ given  as follows:
\be \cd=\{D: \ D= \mathcal A(x), \ x \in X\}. \label{SDP-2}\end{equation}
Evidently, the set $\cd$ is closed,  convex, and the corresponding set of all its zeros (see (\ref{Formula2}))  can be written in the form
\begin{equation} T_0(\cd) =\{t \in T:\ t^\top {\cal A}(x)t=0\; \forall x \in X\}.\label{Tnormalized}\end{equation}

 Suppose that   the    set ${T_0}(\cd)$  is not  empty.
Consider a finite  non-empty subset of this set
\be {\cal V}:=\{{t}(i)\in {T_0}(\cd), \; i \in {\cal I}\},\; 0<|{\cal I}|<\infty.\label{N0**}\ee

\vspace{1mm}

The following theorem permits us to characterize the set (\ref{setX}) in terms  of zeros of the corresponding set $\cd$.}

\begin{theorem}\label{L25-02-1}  Given the closed convex  set  $X$ defined in (\ref{setX}),  the corresponding   subset $\cd$ of the cone $\cop$ defined in (\ref{SDP-2}),  and  any subset  (\ref{N0**}) of the  set $T_0(\cd)$  defined in (\ref{Tnormalized}), the following equality holds:
$ X={\mathcal X}({\cal V}),$
where
$$ \mathcal X({\cal V}):=\{x \in \mathbb R^n: \
{\mathcal A}(x) t(i)\geq 0 \ \forall i \in {\cal I}; \;\;  t^\top {\mathcal A}(x)t\geq 0\; \forall t \in \Omega({\cal V})\},\qquad $$%\label{calX}\end{equation}
and the set $\Omega (\cal V)$ is constructed by the rules (\ref{0**})-(\ref{set-hat1}) with $V={\cal V}.$
\end{theorem}

{\bf Proof.}  Since, by construction, $t(i) \in {T_0}(\cd)$ $ \forall \,i \in \cal I$, it is easy to show that
 $$e^\top_k{\cal A}(x)t(i)=0\;  \forall  k\in P_+(t(i)); \ e^\top_k{\cal A}(x)t(i)\geq 0\;  \forall   k\in P_0(t(i)),\;\forall i \in {\cal I},  \ \forall x\in X.$$%\label{WW1}\ee
%Here and in what follows, $\{e_k,  \ k \in P\}$ is the standard orthogonal basis of the space $\mathbb R^p.$
It follows   from  these relations  that $X\subset {\mathcal X}({\cal V}).$

Now, consider any $x\in {\mathcal X}({\cal V}).$ By construction,
$${\cal A}(x)t(i)\geq 0\; \forall  i \in {\cal I};\ \  t^\top {\cal A}(x)t\geq 0 \; \forall t \in \Omega({\cal V}).$$
It follows from Theorem \ref{LO-1}  that the relations above imply the inequalities $t^\top {\cal A}(x)t\geq 0$ for all $t \in T$  and, consequently, $x \in X.$  Hence, we have shown that
${\mathcal X}({\cal V})\subset  X.$ The theorem is proved. $\; \blacksquare$

\begin{remark} Theorem \ref{L25-02-1} can be considered as a generalization of Lemma 2 from \cite{KT-SetValued}, where it was proved that for the set $X$ defined in (\ref{setX}), there {\bf exists} a number $\sigma_*>0$ such that
\begin{equation*}\begin{split}  X =& \bigl\{x \in \mathbb R^n:  {\mathcal A}(x) t(i)\geq 0\  \forall
i \in {\cal I};\nonumber\\
&\qquad \qquad t^\top {\mathcal A}(x)t\geq 0\; \forall  t \in \{t \in T:\rho(t, {\rm conv} \{t(i),i\in {\cal I}\})\geq \sigma_*\}\bigr\},\end{split}\end{equation*}
 but  no rules  have been given to find this number.
Theorem \ref{L25-02-1} states that one can easily compute  $\sigma_*$ by the formula
$ \ \sigma_*=\sigma({\cal V})=\min\{t_k(i), k \in P_+(t(i)), \  i \in {\cal I}\}, \ $
where ${\cal V}$ is defined in (\ref{N0**}).
\end{remark}

\begin{remark} It is easy to see that  all  the results of this  section remain fair if

$\bullet$  replace the condition $x\in \mathbb R^n$  by the condition $x\in X_*$ with a convex closed set $X_*\subset \mathbb R^n,$

$\bullet$ given a set of indices  $\cal V$  in the form  (\ref{N0**}), replace the number  $\sigma ( {\cal V})$ defined in (\ref{o-1})\ by the number
$ \ \sigma_*( {\cal V}):=\min\{t_k(i),\;  k \in P_+(t(i)),\ i \in  {\cal I}_* \}\geq \sigma ( {\cal V}),\ $
where
 $\{t(i),i \in  {\cal I}_* \}$, $ {\cal I}_* \subset  {\cal I},$ is  the set of vertices of the set ${\rm conv} {\cal V}.$
\end{remark}

  The statement of Theorem \ref{L25-02-1} is  useful for {study of}  linear copositive problems since now we can obtain "good" $\ $ description of the set $X$ and
 regularize this problem.

 In fact, suppose that in problem (\ref{cop}), the matrix function ${\cal A}(x)$ is linear and its  constraints do not satisfy the  Slater condition. Consider  the defined in (\ref{SDP-2}) set $\cd$. Let
$Z_{\cd}:=\{ \tau(j), j\in J_{\cd}\}$ be the set of all  minimal zeros of $\cd$.
It follows from Lemma \ref{LOK} that $ 1\leq  |J_{\cd}|<\infty.$

\vspace{2mm}

 Consider the SIP problem
 $$ \min c(x),\;\;
 {\rm s.t.}\;\; {\cal A}(x)\tau(j)\geq 0\ \forall j\in J_{\cd}, \; t^\top {\cal A}(x)t\geq 0 \;\ \forall \, t \in \Omega(Z_{\cd}).$$
This problem can be considered as a regularization of the original linear CoP problem (\ref{cop}) since
it is equivalent to this problem  (see Theorem \ref{L25-02-1}),
has a  finite number of linear inequality constraints ${\cal A}(x)\tau(j)\geq 0, j\in J_{\cd}$,
and
there  exists $\bar x\in X$ such that $t^\top {\cal A}(\bar x)t>0 \;\ \forall \, t \in \Omega(Z_{\cd})$ (see  {Theorem} \ref{L25-2-2}  in section \ref{S4}  or Lemma 3 in \cite{KT-SetValued}).

\vspace{2mm}

Let $\mathbf{F}_{min}$ be the
smallest (by inclusion) face of ${\cal COP}^p$  containing the  set  $\cal D $  defined by  formula (\ref{SDP-2}) in terms of the constraints of problem (\ref{cop}). The face  $\mathbf{F}_{min}$  is  called (see \cite{Waki}) {\it the minimal face } of this optimization problem.

 For the copositive problem (\ref{cop}), we can formulate the following corollary from  Theorem \ref{T-face} proved in the previous section.

\begin{corollary}\label{Cor3} The minimal face of  problem (\ref{cop}) has the form
$$\mathbf{F}_{min}=\{ D\in \cop: e^\top _kD\tau(j)=0  \ \forall  k \in M_\cd(j), \  \forall  j \in J_\cd\},$$
where  the set $\cd$ is defined in (\ref{SDP-2}),
$\{\tau(j), j\in J_\cd\}$ is  the set  of  all minimal zeros of $\cd$,  and  the sets $ M_\cd(j), j\in J_\cd$, are defined by formulas (\ref{M}) with $Q$ substituted by $\cd$:
$$ M_\cd(j):=\{k\in P:    e^{\top}_k\mathcal A(x)\tau(j)=0 \; \forall x \in X\}, \ j\in J_\cd.$$
\end{corollary}

The importance of Corollary \ref{Cor3} lies in the fact that it  describes the way of representing  the   minimal   face of the copositive problem (\ref{cop}) in an  explicit form  using  the minimal zeros of the set $\cd$ . This representation  may be useful for creating new numerical methods based on the minimal cone representations.

\section{Conclusions}\label{S8}

The  results of the paper permit to explicitly describe  faces of the cone of copositive matrices and the respective dual cones  in terms of the  minimal zeros of these faces. The  novelty of the obtained results consists  of  the fact that they give a new perception of the facial structure of the cone ${\cal COP}^p$ which is not  well investigated yet \cite{Dur}.  In his  paper \cite{Pataki} dedicated to  the  characterization of  more simple cases of convex cones (so-called {\it nice} cones), Gabor Pataki wrote: "...Copositive, and completely positive cones
lie at the other end of the spectrum. Though they are very useful in optimization ..., optimizing over them is more difficult. Also, while considerable progress has been made in describing their geometry ..., a complete understanding
(such as a complete description of their facial structure) is probably out of reach."

We expect that the results of this paper  will   help to better understand the facial structure of the cone $\cop$ and describe its faces explicitly.

    Moreover, the  obtained in  the paper representations are useful for  the  study of convex copositive problems. In particular, they permit to
    \begin{itemize}
  \item  create  regularization procedures  based on the face reduction approach;

\item formulate for CoP problems new optimality conditions without  any CQs;

\item develop new strong duality theory for  copositive optimization based on  an  explicit formulation of the  Extended Lagrange Dual Problem;

    \item   develop numerical methods  for solving  CoP problems.

\end{itemize}

\section*{Acknowledgement}

This work was  supported by the state research program "Convergence"(Republic Belarus), Task 1.3.01,     Portuguese funds through CIDMA - Center for Research and Development in Mathematics and Applications, and FCT -  Portuguese Foundation for Science and Technology, within the project\  UIDB/04106/2020.

\end{document}